\documentclass{article}
\usepackage{amsmath}
\usepackage{amsfonts}
\newcommand{\proof}{\textit{Proof}\textrm{. }}
\newcommand{\epr}{\hfill$\diamondsuit$\medskip}
\newtheorem{theorem}{Theorem}
\newtheorem{corollary}[theorem]{Corollary}
\newtheorem{lemma}[theorem]{Lemma}
\newtheorem{proposition}[theorem]{Proposition}
\title{The centre of generic algebras of small PI algebras}
\author{Thiago Castilho de Mello\thanks{Supported by PhD grant from CNPq},  Plamen Koshlukov\thanks{Partially supported by grants from CNPq
  (No. 304003/2011-5), and  from FAPESP  (No. 2010/50347-9)}
\\
Department of Mathematics, IMECC, UNICAMP\\
S\'ergio Buarque de Holanda 651\\
13083-859 Campinas, SP, Brazil\\
e-mail: \texttt{plamen@ime.unicamp.br}, \texttt{tcmello@ime.unicamp.br}}
\date{}
\begin{document}
\maketitle


\begin{abstract}
Verbally prime algebras are important in PI theory. They are well known over a field $K$ of characteristic zero: 0 and $K\langle T\rangle$ (the trivial ones), $M_n(K)$, $M_n(E)$, $M_{ab}(E)$. Here $K\langle T\rangle$ is the free associative algebra with free generators $T$, $E$ is the infinite dimensional Grassmann algebra over $K$, $M_n(K)$ and $M_n(E)$ are the $n\times n$ matrices over $K$ and over $E$, respectively. Moreover $M_{ab}(E)$ are certain subalgebras of $M_{a+b}(E)$, defined below. The generic algebras of these algebras have been studied extensively. Procesi gave a very tight description of the generic algebra of $M_n(K)$. The situation is rather unclear for the remaining nontrivial verbally prime algebras.

In this paper we study the centre of the generic algebra of $M_{11}(E)$ in two generators. We prove that this centre is a direct sum of the field and a nilpotent ideal (of the generic algebra). We describe the centre of this algebra. As a corollary we obtain that this centre contains nonscalar elements thus we answer a question posed by Berele.
\end{abstract}

The verbally prime algebras (also called T-prime) play a crucial role in the theory of the ideals of identities (also called T-ideals) of associative algebras. A T-ideal is called T-prime if it is prime in the class of all T-ideals. Let $K$ be a field and denote by $K\langle T\rangle$ the free associative algebra freely generated by the set $T$ over $K$. If char$K=0$ then the nontrivial T-prime T-ideals are those of the polynomial identities of the following algebras: $M_n(K)$, $M_n(E)$, $M_{ab}(E)$. We denote here by $E$ the infinite dimensional Grassmann (or exterior) algebra over $K$. The algebra $M_{ab}(E)$ is a subalgebra of $M_{a+b}(E)$. It consists of the block matrices having blocks $a\times a$ and $b\times b$ on the main diagonal with entries from $E_0$, and all remaining entries from $E_1$. Here $E_0$ is the centre of $E$ and $E_1$ is the anticommuting part of $E$. In order to be more precise, assume $V$ is a vector space with a basis $e_1$, $e_2$, \dots, and let $E$ be the Grassmann algebra of $V$. Then $E$ has a basis consisting of all elements of the type $e_{i_1}\ldots e_{i_k}$ where $i_1<\cdots<i_k$, $k\ge 0$, and multiplication induced by $e_ie_j=-e_je_i$. Hence $E_0$ is the span of all of the above elements with even $k$ while $E_1$ is the span of those with odd $k$.

The above classification of the T-prime algebras was obtained by Kemer, as a part of the theory that led him to the positive solution of the Specht problem, see \cite{kemerbook} for an account of Kemer's theory.

Although polynomial identities in T-prime algebras have been extensively studied the concrete information is quite scarce. Thus the polynomial identities for $M_n(K)$ are known only for $n\le 2$, see \cite{razmal, razmbook, dral} when $K$ is of characteristic 0, and \cite{kjam2, jcpk}, when $|K|=\infty$, char$K=p>2$. The identities satisfied by the Grassmann algebra $E$ are well known, see \cite{krreg} when char$K=0$, and the references of \cite{gk} for the remaining cases for $K$. The identities of $M_{11}(E)$ were described in characteristic 0 by Popov, \cite{popov}. Recall that the paper \cite{popov} gives a basis of the identities satisfied by $E\otimes E$ but it is well known (see for example \cite{kemerbook}) that the latter algebra satisfies the same identities as $M_{11}(E)$ when char$K=0$. Our knowledge about the identities even of $M_3(K)$ and of $M_2(E)$ is quite limited; it should be noted that no working methods are available in order to describe them.

Let $A$ be a PI algebra and suppose $I=T(A)$ is its T-ideal in $K\langle T\rangle$. The quotient $K\langle T\rangle/I$ is the relatively free algebra, also called the generic algebra of $A$. Thus one may want to study the generic algebras of the T-prime algebras. It is worth mentioning that these generic algebras admit quite natural models as matrices over certain algebras. The generic algebra for $M_n(K)$ is called the generic matrix algebra. It is a fundamental object in Invariant theory, and enjoys very many good properties; one associates the study of the generic matrix algebras with Procesi, see for example \cite{procesi, procesiam} and also \cite{formanekbook}. Concrete models for the generic algebras for $M_{ab}(E)$ and for $M_n(E)$ were described by Berele \cite{bereleca}. Moreover the study of these generic algebras led to descriptions of their trace rings and to many interesting results in Invariant theory, see for example \cite{bereletams, berelepjm, bereleam, bereleja}. The detailed knowledge of the generic algebra for $M_n(K)$ led Procesi \cite{procesiam} to the description of the trace identities of this algebra, a result obtained independently by Razmyslov as well, see  \cite{razmbook}. The Razmyslov and Procesi's theorem states that the trace identities for $M_n(K)$ all follow from the Cayley--Hamilton characteristic polynomial. Razmyslov proved an analogue of this assertion for the algebras $M_{ab}(E)$ as well.

In \cite[Corollary 21]{bereleca} it was proved that the centre of the generic algebra of $M_{ab}(E)$ is a direct sum of the base field and a nilpotent ideal of the centre. Moreover the author of \cite{bereleca} asked whether the centre of that generic algebra contains any non-scalar elements.

In this paper we describe completely the centre of the generic algebra in two generators of $M_{11}(E)$. It follows from our description that it is a direct sum of the field and a nilpotent ideal (of the generic algebra). Moreover we obtain a detailed information about that nilpotent ideal. As a corollary we show that there are very many non-scalar elements in the centre.

By using this description of the centre we were able to obtain, in characteristic 0, a basis of the polynomial identities satisfied by the generic algebra of $M_{11}(E)$ in two generators. Clearly these differ significantly from the identities of $M_{11}(E)$. This last result requires quite a lot of work, and will be published in a forthcoming paper.

\section{Preliminaries}
We fix an infinite field $K$ of characteristic different from 2. All algebras and vector spaces we consider will  be over $K$. We denote by $K\langle
T\rangle$ the free (unitary) associative algebra freely generated over $K$ by the infinite countable set $T=\{t_1,t_2,\dots\}$. One may conveniently view
$K\langle T\rangle$ as the algebra of polynomials in the non-commuting variables $T$. If $T_k$ is a finite set with $k$ elements, say
$T_k=\{t_1,\dots,t_k\}$ then the free algebra in $k$ generators is denoted by $K\langle T_k\rangle$. The polynomial $f(t_1,\ldots,t_n)\in K\langle
T\rangle$ is a polynomial identity for the algebra $A$ if $f(a_1,\ldots, a_n)=0$ for all $a_i\in A$. The set of all polynomial identities satisfied by $A$
is denoted by $T(A)$, it is its T-ideal. Here we suppose $T(A)\subseteq K\langle T\rangle$. Set $T_k(A) = T(A)\cap K\langle T_k\rangle$. Furthermore we
denote $U(A)= K\langle T\rangle/T(A)$ and $U_k(A)= K\langle T_k\rangle/T_k(A)$ the relatively free algebras of $A$ of infinite rank and of rank $k$,
respectively. With some abuse of notation we shall use the same letters $t_i$ for the free generators of $K\langle T\rangle$ and for their images under the
canonical projection on $U(A)$; analogously for the rank $k$ case.

The algebra $A$ is 2-graded if $A=A_0\oplus A_1$, a direct sum of vector subspaces such that $A_iA_j\subseteq A_{i+j}$ where the latter sum is taken modulo
2. Such algebras are often called superalgebras. A typical example is the Grassmann algebra $E=E_0\oplus E_1$ as above. We call the elements from $A_0\cup A_1$
homogeneous. When $a\in A_i$ we denote its homogeneous degree $\deg a=i$, $i=0$, 1. If $A$ is 2-graded and moreover $ab-(-1)^{\deg a\deg b} ba=0$ for all
homogeneous $a$ and $b$ then $A$ is called a supercommutative algebra. Clearly the Grassmann algebra is supercommutative. Next we recall the construction
of the free supercommutative algebra, see for example \cite[Lemma 1]{bereleca}. Let $X$ and $Y$ be two sets and form the free associative algebra $K\langle X\cup
Y\rangle$. It is 2-graded assuming the elements of $X$ of degree 0 and those of $Y$ of degree 1. Denote by $I$ the ideal generated by all $ab -(-1)^{\deg
a\deg b}ba$ where $a$, $b$ are homogeneous, and put $K[X;Y]=K\langle X\cup Y\rangle/I$. It is immediate to see that $K[X;Y]\cong K[X]\otimes_K E(Y)$. Here
$K[X]$ is the polynomial algebra in $X$ and $E(Y)$ is the Grassmann algebra of the vector space with basis $Y$. Thus if $Y=\{y_1,y_2,\ldots\}$ then $E(Y)$ will have a basis consisting of the products $y_{i_1}\cdots y_{i_k}$, $i_1<\cdots<i_k$, and multiplication induced by $y_iy_j=-y_jy_i$. Below we also recall the construction of
the generic algebras for the T-prime algebras.

Suppose $X=\{x_{ij}^r\}$, $Y=\{y_{ij}^r\}$ where $1\le i,j\le n$, $r=1$, 2, \dots; observe that we use $r$ as an upper index, not as an exponent. Define
the matrices $A_r=(x_{ij}^r)$, $B_r=(x_{ij}^r+y_{ij}^r)$, $C_r=(z_{ij}^r)$ where $z=x$ whenever $1\le i,j\le a$ or $a+1\le i,j\le a+b$, and $z=y$ for all
remaining possibilities for $i$ and $j$. Suppose $a+b=n$, and consider the following subalgebras of $M_n(K[X;Y])$. The first is generated by the generic
matrices $A_r$, $K[A_r\mid r\ge 1]$. It is well known it is isomorphic to the relatively free (or universal) algebra $U(M_n(K))$ of $M_n(K)$. In
\cite[Theorem 2]{bereleca} it was shown that $U(M_n(E))\cong K[B_r\mid r\ge 1]$, and that $U(M_{ab}(E))\cong K[C_r\mid r\ge 1]$. Moreover the relatively
free algebras of finite rank $k$, denoted by $U_k$, can be obtained by letting $r=1$, \dots, $k$, that is by taking the first $k$ matrices.

We recall another fact from \cite{bereleca} that we shall exploit. It was shown in \cite[Theorem 20]{bereleca} that if $f$ is a central polynomial for $M_{ab}(E)$, without constant term, then for some $m$ the polynomial $f^m$ is an identity for $M_{ab}(E)$. It follows that the centre of $U_k(M_{ab}(E))$ must be a direct sum of $K$ and a nilpotent ideal of the centre, see \cite[Corollary 21]{bereleca}.

We shall need information about the polynomial identities of $M_{11}(E)$. These were described by Popov in characteristic 0, see the main theorem of \cite{popov}. As we mentioned above, in \cite{popov} it was proved that the T-ideal of $E\otimes E$ is generated by the two polynomials
\begin{equation}
\label{basispim11}
[[t_1,t_2]^2,t_1], \qquad [[t_1,t_2],[t_3,t_4],t_5]
\end{equation}
where $[a,b]=ab-ba$ is the usual commutator. We consider the commutators left normed that is $[a,b,c] = [[a,b],c]$, and so on in higher degree.

The algebra $K\langle T\rangle$ is multigraded by the degree of its monomials in each variable.
We work with the infinite field $K$ therefore every T-ideal is generated by its multihomogeneous elements, see for example \cite[Section 4.2]{drenskybook}. Thus from now on we shall work with multihomogeneous polynomials only.

The algebra $E\otimes E$ is PI equivalent to $M_{11}(E)$ in characteristic 0, so the polynomials (\ref{basispim11}) generate the T-ideal of $M_{11}(E)$ as well. This is a result due to Kemer, see \cite{kemerbook}. Kemer proved that the tensor product of two T-prime algebras (in characteristic 0) is PI equivalent to a T-prime algebra, and described precisely these PI equivalences.  Note that if char$K=p>2$ then the algebras $M_{11}(E)$ and $E\otimes E$ are not PI equivalent, see for example \cite{afk}, or \cite{smapk}. While the former paper proved directly the non-equivalence the latter proved it by computing the GK dimensions of the corresponding relatively free algebras (these turn out to be different).

\section{The free supercommutative algebra}

In this section we denote by $F=U_2(M_{11}(E))$ the relatively free algebra of rank 2 for $M_{11}(E)$. As mentioned before we have that $F=K[C_1, C_2]$. From now on
we consider the free supercommutative algebra $K[X;Y]$ where $X=\{x_1,x_2, x_1',x_2'\}$, $Y=\{y_1,y_2,y_1',y_2'\}$, and set $C_1=\begin{pmatrix}
x_1&y_1\\ y_1'& x_1'\end{pmatrix}$, $C_2=\begin{pmatrix} x_2&y_2\\ y_2'& x_2'\end{pmatrix}$. Clearly the algebra $F$ satisfies all identities of
$M_{11}(E)$.

The algebra $K[X;Y]$ is graded by the integers, taking into account the degree with respect to the variables in $Y$ only: $K[X;Y] = \oplus_{n\in\mathbb{Z}}
K[X;Y]^{(n)}$. Here $K[X;Y]^{(n)}$ is the span of all monomials of degree $n$ in the variables from $Y$. It is immediate that the $n$-th homogeneous component is zero
unless $0\le n\le 4$. The canonical 2-grading on $K[X;Y]$ and the $\mathbb{Z}$-grading just defined are related as follows:
\begin{eqnarray*}
K[X;Y]_0 &=& K[X;Y]^{(0)} + K[X;Y]^{(2)} + K[X;Y]^{(4)}; \\
K[X;Y]_1 &=& K[X;Y]^{(1)} + K[X;Y]^{(3)}.
\end{eqnarray*}

The next facts are quite obvious; we collect them in a lemma for further reference.

\begin{lemma}
\label{obvious} Consider the polynomial algebra $K[X]\subseteq K[X;Y]$.

1. For every $n$, $0\le n\le 4$, $K[X;Y]^{(n)}$ is a free module over $K[X]$, with a basis $B_n$, where $B_0=\{1\}$, and
\begin{eqnarray*}
&&B_1=\{y_1,y_2,y_1', y_2'\}, \quad
B_2=\{y_1y_2, y_1y_1', y_1y_2', y_2y_1', y_2y_2', y_1'y_2'\}\\
&&B_3=\{y_1y_2y_1', y_1y_2y_2', y_1y_1'y_2', y_2y_1'y_2'\},  \quad B_4=\{y_1y_2y_1'y_2'\}.
\end{eqnarray*}

2. The free supercommutative algebra $K[X;Y]$ is a free module over $K[X]$ with a basis $B=B_0\cup B_1\cup B_2\cup B_3\cup B_4$.

3. Every ideal of $K[X;Y]$ is a $K[X]$-submodule of $K[X;Y]$.
\end{lemma}

In effect one may extend the scalars as follows. Let $K(X)$ be the field of fractions of $K[X]$, and consider the Grassmann algebra $\overline{E(Y)}$ on $Y$ over $K(X)$, that is $\overline{E(Y)} = E(Y)\otimes_K K(X)$.

\begin{lemma}
\label{nozerodiv}
The matrices $C_1$ and $C_2$ are not zero divisors in $F$.
\end{lemma}

\proof
Suppose $A=\begin{pmatrix} a&b\\ c&d\end{pmatrix}\in F$ is such that $C_1A=0$. Then one obtains
\[
x_1a+y_1c=0, \quad y_1'a+x_1'c=0, \quad x_1b+y_1d=0, \quad y_1'b+x_1'd=0.
\]
It follows from the second equation that $x_1'c=-y_1'a$. Now multiply the first equation by $x_1'$ and substitute $x_1'c$ by its equal and get $a(x_1x_1'-y_1y_1')=0$. Extending the scalars as above, and noting that the set $B$ from Lemma~\ref{obvious} is a basis of the vector space $\overline{E(Y)}$ over $K(X)$, we immediately obtain $a=0$ and consequently $c=0$. In the same manner but using the last two equations we obtain $b=d=0$. In this way $C_1$ is not a left zero divisor. Analogously one shows it is not a right zero divisor, and the same for $C_2$.
\epr

We define an automorphism $'$ on $K[X;Y]$ by setting $x_i\mapsto x_i'$, $x_i'\mapsto x_i$, $y_i\mapsto y_i'$ and $y_i'\mapsto y_i$. Thus the automorphism $'$ is of order two.

We also define the following polynomials in $K[X]$:
\begin{eqnarray*}
&&q_n(x_1,x_1') = \sum_{i=0}^nx_1^i x_1'^{n-i}; \quad
Q_n(x_2,x_2') = q_n(x_2,x_2');\\
&&r_n(x_1, x_1') = \sum_{i=0}^{n-1} (n-i)x_1^{n-1-i}x_1'^i; \quad
R_n(x_2,x_2') = r_n(x_2,x_2');\\
&&s_n(x_1,x_1') = r_n(x_1', x_1); \quad
S_n(x_2,x_2') =  s_n(x_2,x_2').
\end{eqnarray*}

\begin{lemma}
\label{relationsqrs}
The following relations hold among the above polynomials.
\begin{eqnarray*}
&&r_n=q_{n-1}+x_1r_{n-1}; \quad s_n=q_{n-1} + x_1's_{n-1}; \quad s_n+r_n = (n+1) q_{n-1}; \\
&&q_n=x_1^n + x_1' q_{n-1} = x_1'{}^n + x_1q_{n-1};\quad (x_1'-x_1)q_{n-1} = x_1'{}^n - x_1^n;\\
&& x_1^nx_1'{}^m - x_1^mx_1'{}^n = (x_1'-x_1)(q_nq_{m-1} - q_mq_{n-1}).
\end{eqnarray*}
\end{lemma}

\proof
The proof consists of an easy induction.
\epr

It is immediate to check, once again by induction,  that for every $m$ and $n$,
\begin{eqnarray*}
C_1^n &=&\begin{pmatrix}
x_1^n+y_1y_1'r_{n-1}& y_1q_{n-1}\\
y_1'q_{n-1}& x_1'{}^n-y_1y_1's_{n-1}
\end{pmatrix}; \\
C_2^m &=& \begin{pmatrix}
x_2^m + y_2 y_2' R_{m-1} & y_2 Q_{m-1}\\
y_2' Q_{m-1} & x_2'{}^m - y_2y_2' S_{m-1}
\end{pmatrix}.
\end{eqnarray*}
Note that the $y_i$ and the $y_i'$ anticommute and this produces the minus signs at the $(2,2)$-entries of the above matrices.

Therefore for the product $C_1^nC_2^m$ we have
\[
C_1^nC_2^m = \begin{pmatrix}
x_1^nx_2^m+a+d & y_1x_2'{}^mq_{n-1} + y_2x_1^n Q_{m-1} + c\\
y_1'x_2^m q_{n-1} + y_2'x_2'{}^nQ_{m-1} + c' & x_1'{}^n x_2'{}^m+ a'+d'
\end{pmatrix}
\]
where $a$, $a'\in K[X;Y]^{(2)}$, $d$, $d'\in K[X;Y]^{(4)}$, and $c$, $c'\in K[X;Y]^{(3)}$.

As the elements of the algebra $F$ are linear combinations of products of the above type we obtain immediately the proof of the following lemma.

\begin{lemma}
\label{autoorder2} Let $A=(a_{ij})\in F$, then $a_{22}=a_{11}'$ and $a_{21}=a_{12}'$.
\end{lemma}

We shall need the following elements in order to describe the centre of $F$.
\begin{eqnarray*}
h_1&=& y_1y_2y_1'y_2';\\
h_2 &=& y_1y_2 (y_1'(x_2'-x_2) - y_2'(x_1'-x_1)); \\
h_3 &=& y_1'y_2' (y_1(x_2'-x_2) - y_2(x_1'-x_1)); \\
h_4 &=& (y_1'(x_2'-x_2) - y_2'(x_1'-x_1))
(y_1(x_2'-x_2) - y_2(x_1'-x_1)).
\end{eqnarray*}

Once again using the fact that $y_1$, $y_1'$, $y_2$, $y_2'$ anticommute we obtain that the elements $h_1$, $h_2$, $h_3$, $h_4$ satisfy the following relations in $K[X;Y]$.
\begin{eqnarray*}
&&h_1y_1 = h_1y_1' = h_1y_2 = h_1y_2' = 0; \quad h_2y_1 = h_2y_2 = h_3y_1' = h_3y_2' = 0;\\
&&h_2y_1' = h_3y_1 = (x_1'-x_1)h_1; \quad h_2y_2' = h_3y_2 = (x_2'-x_2) h_1; \\
&&h_4y_1 = (x_1'-x_1)h_2; \quad h_4y_2 = (x_2'-x_2) h_2; \\
&&h_4y_1' = -(x_1'-x_1) h_3; \quad h_4y_2' = -(x_2'-x_2)h_3.
\end{eqnarray*}

\section{The centre of $F_2(M_{11}(E))$}

Take a matrix $A=\begin{pmatrix} a&b\\ c&d\end{pmatrix}\in F\cong K[C_1,C_2]$. Assume it is central in $F$, that is  $[A,C_1]=[A, C_2]=0$. Then
\begin{eqnarray}
\label{matrixentries}
&&by_1'+cy_1=0; \quad by_2'+cy_2=0; \notag\\
&&(a-d)y_1+(x_1'-x_1)b=0; \quad
(a-d)y_2+(x_2'-x_2)b=0;\\
&&(a-d)y_1'+(x_1'-x_1)c=0;
\quad
(a-d)y_2'+(x_2'-x_2)c=0.\notag
\end{eqnarray}
We multiply the third equation by $(x_2'-x_2)$, the fourth by $(x_1'-x_1)$ and subtract, obtaining $(d-a)((x_2'-x_2)y_1 - (x_1'-x_1)y_2) = 0$. Similarly
from the last two equations we get $(d-a)((x_2'-x_2)y_1' - (x_1'-x_1)y_2') = 0$. Therefore
\[
d-a\in Ann ((x_2'-x_2)y_1 - (x_1'-x_1)y_2)\cap
Ann ((x_2'-x_2)y_1' - (x_1'-x_1)y_2').
\]

Denote by $J$ the intersection of the two annihilators in the right hand side above.
Then $J$ is an ideal of $K[X;Y]$ hence $J$ is a $K[X]$-submodule as well.

\begin{proposition}
\label{basisJ}
The $K[X]$-module $J$ is spanned by $\{h_1,h_2,h_3,h_4\}$.
\end{proposition}

\proof It is immediate that $J$ contains $h_1$, $h_2$, $h_3$, $h_4$.  We shall prove that $J$ is contained in the $K[X]$-module spanned by
$\{h_1,h_2,h_3,h_4\}$. Let $f\in J$ and write $f=f_0+f_1+f_2+f_3+f_4$ where $f_i\in K[X;Y]^{(i)}$. First note that both $((x_2'-x_2)y_1 - (x_1'-x_1)y_2)$ and
$((x_2'-x_2)y_1' - (x_1'-x_1)y_2')$ lie in $K[X;Y]^{(1)}$ in the $\mathbb{Z}$-grading defined above. Thus $f$ annihilates the latter two polynomials if and
only if every $f_i$ does. Therefore we may and shall assume $f$ is homogeneous in the $\mathbb{Z}$-grading.

Suppose first $f\in K[X;Y]^{(0)}\cong K[X]$. As $f\in J$ then it follows easily that $f=0$.

Let $f\in K[X;Y]^{(1)}\cap J$, then $f=\alpha_1 y_1 + \alpha_2 y_2 + \alpha_3 y_1' + \alpha_4 y_2'$ for some $\alpha_i\in K[X]$. But $f(y_1(x_2'-x_2) -
y_2(x_1'-x_1))=0$, hence
\[
(x_1'-x_1)(-\alpha_1y_1y_2 + \alpha_3 y_2y_1' +\alpha_4 y_2y_2') - (x_2-x_2')(-\alpha_2 y_1y_2 - \alpha_3y_1y_1' - \alpha_4 y_1y_2'))=0.
\]
As $B_2$ is a $K[X]$-basis for the free module $K[X;Y]^{(2)}$ it follows $\alpha_3=\alpha_4=0$. In the same way, using the fact that $f$ annihilates
$y_1'(x_2'-x_2) - y_2'(x_1'-x_1)$, we get $\alpha_1=\alpha_2=0$. This implies $J\cap K[X;Y]^{(1)}=0$.

Let $f= \alpha_1 y_1y_2 + \alpha_2y_1y_1' + \alpha_3 y_1y_2' + \alpha_4 y_2y_1' + \alpha_5 y_2y_2' + \alpha_6 y_1'y_2'\in K[X;Y]^{(2)}\cap J$. As above,
from $f(y_1(x_2'-x_2) - y_2(x_1'-x_1))=0$ we obtain, in $K[X;Y]^{(3)}$, a linear combination of $y_1y_2y_1'$, $y_1y_2y_2'$, $y_1y_1'y_2'$, $y_2y_1'y_2'$,
with coefficients respectively
\[
\alpha_2(x_1'-x_1) + \alpha_4(x_2'-x_2);
\alpha_3(x_1'-x_1) + \alpha_5(x_2'-x_2);
\alpha_6 (x_2'-x_2); -\alpha_6 (x_1'-x_1).
\]
Thus $\alpha_6=0$, $\alpha_2(x_1'-x_1) + \alpha_4(x_2'-x_2)=0$, and $\alpha_3(x_1'-x_1) + \alpha_5(x_2'-x_2)=0$.

Analogously, from  $f(y_1'(x_2'-x_2) - y_2'(x_1'-x_1))=0$ we obtain that $\alpha_1=0$,  $\alpha_2(x_1'-x_1) + \alpha_3(x_2'-x_2)=0$, and $\alpha_4(x_1'-x_1) + \alpha_5(x_2'-x_2)=0$.

Now the $\alpha_i$ are polynomials in $K[X]$. Consider the field of fractions $K(X)$ of this polynomial ring, and resolve the corresponding linear system of four equations in $K(X)$. One obtains that the solution depends on one parameter $\beta\in K(X)$:
\[
\alpha_2=\beta; \quad \alpha_3 = -\frac{x_1'-x_1}{x_2'-x_2}\beta; \quad \alpha_4 = -\frac{x_1'-x_1}{x_2'-x_2}\beta; \quad \alpha_5=
\frac{(x_1'-x_1)^2}{(x_2'- x_2)^2}\beta,
\]
and of course $\alpha_1=\alpha_6=0$. Since we are looking for a solution of the system in $K[X]$ then one must have $\beta = (x_2'-x_2)^2\alpha$ for some $\alpha\in K[X]$, and the solution in $K[X]$ will be
\begin{eqnarray*}
&& \alpha_2= (x_2'-x_2)^2 \alpha; \quad
\alpha_3 = -(x_1'-x_1)(x_2'-x_2) \alpha; \\
&&\alpha_4 = -(x_1'-x_1)(x_2'-x_2) \alpha; \quad
\alpha_5 = (x_1'-x_1)^2 \alpha.
\end{eqnarray*}
When we substitute these in the expression of $f$ we get $f=\alpha h_4$.

Let $f=\alpha_1 y_1y_2y_1' + \alpha_2 y_1y_2y_2' + \alpha_3 y_1y_1'y_2' + \alpha_4 y_2y_1'y_2'\in K[X;Y]^{(3)}\cap J$ where $\alpha_i\in K[X]$. Proceeding
as above we get the equalities
\begin{eqnarray*}
(\alpha_4(x_2'-x_2)+\alpha_3 (x_1'-x_1)) y_1y_1'y_2y_2' &=&0\\
(\alpha_2(x_2'-x_2)+\alpha_1 (x_1'-x_1)) y_1y_1'y_2y_2' &=&0
\end{eqnarray*}
therefore $\alpha_2(x_2'-x_2)+\alpha_1 (x_1'-x_1) = 0$ and $\alpha_4(x_2'-x_2)+\alpha_3 (x_1'-x_1) =0$. As above we work first in $K(X)$ and then go back to $K[X]$. We find that the solution in $K[X]$ is
\[
\alpha_1 = -(x_2'-x_2)\alpha; \quad \alpha_2 = (x_1'-x_1)\alpha; \quad \alpha_3 = -(x_2'-x_2) \beta; \quad \alpha_4= (x_1'-x_1)\beta
\]
where $\alpha$, $\beta\in K[X]$. Substituting these values in $f$ we obtain $f=\alpha h_2+ \beta h_3$, and this case is dealt with.

Finally it is immediate to see that $K[X;Y]^{(4)} = K[X]\cdot h_1$, and thus we conclude the proof. \epr

\begin{corollary}
\label{conditionsonelements}
Let the matrix $A$ be as at the beginning of this section. Then $A$ commutes with $C_1$ and $C_2$ if and only if $b=f_4 h_2$, $c = -f_4 h_3$, $d = a+f_1h_1+f_4h_4$ for some $f_1$, $f_4\in K[X]$.
\end{corollary}

\proof We already saw that $d-a\in J$. Thus $d-a = f_1h_1+f_2h_2+f_3h_3+f_4h_4$, $f_i\in K[X]$. By means of homogeneity we get $d-a\in K[X;Y]_{0}$, it
follows that $d-a=f_1h_1+f_4h_4$. Now substituting $d-a$ in the system (\ref{matrixentries}) we have $b=f_4 h_2$ and $c= -f_4 h_3$. On the other hand if
$a$, $b$, $c$, $d$ satisfy the conditions of the statement it is immediate that $[A,C_1] = [A,C_2] = 0$. \epr

\noindent
\textbf{Remark}
1. We just proved that $[A,C_1] = [A,C_2] = 0$ if and only if
\[
A= aI + f_1\begin{pmatrix} 0&0\\ 0&h_1\end{pmatrix} + f_4\begin{pmatrix} 0&h_2\\ -h_3& h_4\end{pmatrix}.
\]

2. Therefore if $A$ is central in $F$ then $a_{12}$, $a_{21}\in K[X;Y]^{(3)}$.

\bigskip

An element $a\in F$ will be called \textsl{strongly central} if it is central, and moreover, for every $b\in F$ the element $ab$ is central in $F$ (thus $ba=ab$ will be strongly central as well).

Let us fix the following matrices in $F$:
\[
A_0=\begin{pmatrix} h_1&0\\ 0&h_1\end{pmatrix}; \quad
A_1=\begin{pmatrix} 0&0\\ 0&h_1\end{pmatrix}; \quad
A_2=\begin{pmatrix} 0&h_2\\ -h_3&h_4\end{pmatrix}; \quad
A_3=\begin{pmatrix} h_4&0\\ 0&h_4\end{pmatrix}.
\]

\begin{lemma}
\label{stronglycentral} Let $a=\alpha_0A_0 + \alpha_1A_1 + \alpha_2A_2 + \alpha_3A_3$, for some $\alpha_i\in K[X]$. If $a\in F$ then $a$ is strongly
central.
\end{lemma}

\proof
The matrices $A_0$ and $A_1$ are clearly strongly central. Also $A_2$ and $A_3$ are central. One computes
\[
A_2C_i = (x_i'-x_i)\begin{pmatrix} h_1&0\\ 0&-h_1 \end{pmatrix} + x_i'\begin{pmatrix} 0&h_2\\ -h_3& h_4\end{pmatrix}.
\]
Hence $A_2C_i$ is a linear combination (over $K[X]$) of $A_0$, $A_1$ and $A_2$. Iterating we will have $A_2C_{i_1}C_{i_2}\dots C_{i_r}$ is central for
$i_j=1,2$ and $r=1,2\dots$, and $A_2$ is strongly central. One checks in a similar manner that

\[
A_3C_i = x_i\begin{pmatrix} h_4&0\\ 0&h_4\end{pmatrix} + (x_i'-x_i)\begin{pmatrix} 0&h_2\\ -h_3&h_4\end{pmatrix}.
\]
That is $A_3C_i$ is a combination of $A_2$ and $A_3$ and iterating as above we show that $A_3$
is strongly central. \epr

\begin{lemma}
\label{commutator}
Let $f(t_1,t_2)=[t_1,t_2,t_{i_3}, \ldots, t_{i_k}]$ be a left normed commutator, $i_j=1$, 2. Suppose that $\deg_{t_1} f=n$, $\deg_{t_2} f=m$, $n+m=k$. Then for every $n$ and $m$ one has $f(C_1,C_2) = (x_1'-x_1)^{n-1} (x_2'-x_2)^{m-1} A(k)$ where
\begin{eqnarray*}
A(k) &=& \begin{pmatrix} F(k) & y_1(x_2'-x_2) - y_2 (x_1'-x_1)\\
(-1)^k (y_2'(x_1'-x_1) - y_1'(x_2'-x_2))& F(k) \end{pmatrix},\\
F(k) &=& \frac{(y_1(x_2'-x_2) - y_2(x_1'-x_1)) y_i' + (-1)^k y_i (y_2'(x_1'-x_1) - y_1'(x_2'-x_2))}{x_i'-x_i} .
\end{eqnarray*}
In the last expression we use the shorthand $i$ for $i_k$.
\end{lemma}

\proof
The proof consists of an induction on $k$. The base of the induction is $k=2$; then $F(2) = y_1y_2'+y_1'y_2$. If $f=f_k$ is the commutator of the statement then one computes $[f_k,C_i]$ directly by induction.
\epr

\noindent
\textbf{Remark}
It follows from the above lemma that if $u$ is a left normed commutator in $t_1$ and $t_2$ then $u(C_1,C_2)$ does not depend on the order of the variables starting with the third and up to the last but one. In other words any permutation of the variables in $u$ that preserves the first two and the last one, leaves $u$ invariant.

\begin{lemma}
\label{prod2comm} Let $u_1(t_1,t_2)$ and $u_2(t_1,t_2)$ be two left-normed commutators of degrees at least two in $F$, and denote by  $u=u_1u_2$ their product. Then $u$ is strongly central in $F$.
\end{lemma}

\proof
Suppose $\deg_{t_1} u_j=n_j$ and $\deg_{t_2} u_j = m_j$, $j=1$, 2. Using the notation of  Lemma~\ref{commutator} we have
\[
u_1(C_1,C_2)u_2(C_1,C_2) = (x_1'-x_1)^{n_1+n_2-2} (x_2'-x_2)^{m_1+m_2-2} A(k_1) A(k_2)
\]
where $\deg u_j=k_j$. But it is immediate to see that $F(k_1)F(k_2) = \alpha h_1$, and also
\begin{eqnarray*}
&&F(k_i)(y_1(x_2'-x_2)-y_2(x_1'-x_1)) = -(-1)^{k_i} h_2, \\
&&F(k_i)(y_2'(x_1'-x_1)-y_1'(x_2'-x_2)) = h_3
\end{eqnarray*}
where $\alpha\in K[X]$. Then
\[
A(k_1)A(k_2) = F(k_1)F(k_2) I +\begin{pmatrix} (-1)^{k_2} h_4 & -((-1)^{k_1} + (-1)^{k_2}) h_2\\
((-1)^{k_1} + (-1)^{k_2}) h_3 & -(-1)^{k_1} h_4
\end{pmatrix}.
\]
Therefore $u_1(C_1,C_2)u_2(C_1,C_2) = \alpha A_0 + (-1)^{k_2} A_3 - ((-1)^{k_1} + (-1)^{k_2}) A_2$ is strongly central in $F$.
\epr

\noindent
\textbf{Remark} Let $u_1$ and $u_2$ be two commutators, $\deg u_1\equiv\deg u_2\pmod{2}$. Then $u=u_1u_2$ is central element in $F$ but it is not a scalar multiple of $I$. In particular
\[
[C_1,C_2]^2 = \begin{pmatrix} -2h_1 + h_4& -2h_2\\ 2h_3 & -2h_1 - h_4\end{pmatrix}
\]
is central in $F$ but is not a scalar. This answers Berele's question from \cite{bereleca} for the case of $M_{11}(E)$ and two generators. In this same case we shall give below the precise answer to Berele's question.

\medskip

We consider unitary algebras. Let $L(T)$ be the free Lie algebra freely generated by $T$; suppose further $L(T)\subseteq K\langle T\rangle$. That is we consider the vector space $K\langle T\rangle$ with the commutator operation $[a,b]=ab-ba$, and take $L(T)$ as the Lie subalgebra generated by $T$. Choose an ordered basis of $L(T)$ such that the variables from $T$ precede the longer commutators. As $K\langle T\rangle$ is the universal enveloping algebra of $L(T)$ one has that a basis of $K\langle T\rangle$ consists of 1 and all products $t_{i_1}^{n_1}\cdots t_{i_k}^{n_k} u_{j_1}\cdots u_{j_m}$ where $i_1<\cdots<i_k$, and the $u_{j_i}$ are commutators of degree at least two. Clearly all this holds for $K\langle t_1,t_2\rangle$ and for its homomorphic image $F=K[C_1,C_2]$. Therefore every element of $F$ is a linear combination of products of the type $C_1^n C_2^m u_1\cdots u_r$ where the $u_i$ are commutators. Moreover we can assume all commutators left normed, and of the type $[C_1,C_2, \ldots]$.

\begin{proposition}
\label{centralone}
Let $f(C_1,C_2) = C_1^n C_2^m u_1^{k_1}\cdots u_r^{k_r} \in F$ where the $u_i$ are left normed commutators, $\deg u_i\ge 2$. The element $f$ is central in $F$ if and only if $k_1+\cdots+k_r\ge 2$, or else $m=n=k_1=\cdots=k_r=0$.
\end{proposition}

\proof It follows from Lemma~\ref{prod2comm} that the product of two commutators is strongly central. Thus if $k_1+\cdots+k_r\ge 2$ then $f(C_1,C_2)$ is
central. It remains to prove that $C_1^n C_2^m$ and $C_1^nC_2^m u$ are not central where $u$ is a left normed commutator. (In the former we assume
$n+m>0$.) Take first $C_1^nC_2^m$. By the form of the product (computed just before Lemma~\ref{autoorder2}) it follows that the $(1,2)$ entry of
$C_1^nC_2^m$ is $\alpha y_1+\beta y_2+\gamma$ where $(\alpha,\beta)\ne (0,0)$ and $\gamma\in K[X;Y]^{(3)}$. Now by the remark following
Corollary~\ref{conditionsonelements} we have that $C_1^nC_2^m$ cannot be central as long as $n+m>0$.

One proceeds in a similar manner when $f=C_1^nC_2^m u$ where $u$ is a left normed commutator. The $(1,2)$ entry of $f$ will be $x_1^nx_2^m(y_1(x_2'-x_2) -
y_2(x_1'-x_1))+b$ for some $b\in K[X;Y]^{(3)}$. Once again the remark mentioned above yields that $f$ cannot be central. \epr

\begin{proposition}
\label{equivalent} Let $u_1$, \dots, $u_r$ be left normed commutators, $\deg_{t_1}u_j=n$, $\deg_{t_2} u_j=m$ for all $j$. Suppose $u_j=[t_1,t_2,
t_{j_3},\ldots, t_{j_k}]$, $j=1$, \dots, $r$. Put $f(C_1,C_2) = \sum_j \alpha_j u_j(C_1,C_2)\in F$ where $\alpha_j\in K[X]$. Then the following three
conditions are equivalent.
\begin{itemize}
\item[(1)]
The element $f(C_1,C_2)$ is strongly central in $F$.
\item[(2)]
The element $f(C_1,C_2)$ is central in $F$.
\item[(3)]
The sum $\alpha_1+\cdots+\alpha_r=0$ in $K[X;Y]$.
\end{itemize}
\end{proposition}

\proof
Clearly $(1)$ implies $(2)$. We prove now that $(2)$ implies $(3)$. Suppose $f(C_1,C_2)$ is central in $F$. By Lemma~\ref{commutator} the $(1,2)$ entry of every commutator $u_i$ equals $(x_1'-x_1)^{n-1}(x_2'-x_2)^{m-1} (y_1(x_2'-x_2) - y_2(x_1'-x_1))$. Hence the $(1,2)$ entry of $f(C_1,C_2)$ equals
\[
\sum\beta_i (y_1(x_2'-x_2) - y_2(x_1'-x_1)) = (y_1(x_2'-x_2) - y_2(x_1'-x_1))\sum\beta_i
\]
where $\beta_i = (x_1'-x_1)^{n-1}(x_2'-x_2)^{m-1}\alpha_i$. Thus if $\sum\alpha_i\ne 0$ then the $(1,2)$ entry of $f(C_1,C_2)$ is a non-zero multiple of
$y_1(x_2'-x_2) - y_2(x_1'-x_1)$ by some element of $K[X;Y]_0$ and cannot belong to $K[X;Y]^{(3)}$. By the remark following
Corollary~\ref{conditionsonelements}, $f(C_1,C_2)$ cannot be central.

In order to complete the proof we have to prove that $(3)$ implies $(1)$. Suppose $\sum\alpha_i=0$. It was observed in the remark preceding Lemma~\ref{prod2comm} that if $u_i$ and $u_j$ have the same rightmost variable then $u_i(C_1,C_2) = u_j(C_1,C_2)$. Thus we divide the commutators $u_j$ into two types according to their rightmost variable. Clearly if all of them end with say $t_1$ then $\sum \alpha_j u_j(C_1,C_2) = u_1(C_1,C_2)\sum\alpha_j=0$. Hence suppose $u_1$ ends with $t_1$ while $u_2$ ends with $t_2$. Write $\sum\alpha_ju_j = \beta_1 u_1+ \beta_2 u_2$ where $\beta_q$ is the sum of all $\alpha_j$ such that $u_j$ ends with $t_q$, $q=1$, 2. Then $\beta_1+\beta_2=\sum\alpha_j=0$ and it suffices to prove $u_1-u_2$ is strongly central. But
$u_1-u_2=(x_1'-x_1)^{n-1} (x_2'-x_2)^{m-1} (F_1(k) - F_2(k))I$ where we denote by $F_j(k)$ the expression $F(k)$ from Lemma~\ref{commutator} obtained by $u_j$, $j=1$, 2. Clearly $F_1(k)-F_2(k)=0$ if $k$ is even, and $F_1(k)-F_2(k)= -2(x_1'-x_1)^{-1}(x_2'-x_2)^{-1} h_4$ if $k$ is odd. In this way either $u_1-u_2=0$ or $u_1-u_2$ is a multiple of $A_3$. In both cases it is strongly central in $F$.
\epr

\noindent\textbf{Remark }
We observe that in the previous proposition if we suppose $\alpha_j\in K[X;Y]_0$, the statement of the proposition remains valid replacing the condition
(3) by the condition

(3') The sum $\alpha_1+\cdots +\alpha_r\in K[X;Y]^{(4)}$.

\medskip

Let $f(C_1,C_2)\in F=K[C_1,C_2]$. Then $f$ can be written as
\[
f(C_1,C_2) = \sum_{n,m\ge 0} \alpha_{nm} C_1^nC_2^m + \sum_{n_j,m_j} C_1^{n_j}C_2^{m_j} \sum_i \beta_{ij}u_{ij} + g(C_1,C_2).
\]
Here $\alpha_{nm}$, $\beta_{ij}\in K$, $u_{ij}$ are left normed commutators as in Proposition~\ref{equivalent},
and moreover $g(C_1,C_2) = \sum_u \gamma_u C_1^nC_2^m u_1\cdots u_k$ where $u_i$ are left normed commutators. Define $I_{ij}$ as the set of all indices $p$
such that $\deg_{t_1} u_{pj} = r_{ij}$, $\deg_{t_2} u_{pj} = s_{ij}$ for some integers $r_{ij}$ and $s_{ij}$.

\begin{theorem}
\label{descriptioncentre}
Using the notation above, $f(C_1,C_2)$ is central in $F$ if and only if $\alpha_{nm}=0$ for all $n$ and $m$ such that $n+m\ge 1$, and moreover, for  every $i$ and $j$ the equalities $\sum_p \beta_{pj} = 0$ hold where $p\in I_{ij}$.

Furthermore $f(C_1, C_2)$ is strongly central if and only if it is central and $\alpha_{00}=0$.
\end{theorem}

\proof
We already proved that $f(C_1,C_2)$ is central provided that $\alpha_{nm}=0$ when $n+m\ge 1$ and all sums $\sum_{p\in I_{ij}} \beta_{pj}=0$. Such an element is strongly central if and only if $\alpha_{00}=0$. We shall prove the converse. Clearly $g(C_1,C_2)$ is strongly central and $\alpha_{00}I$ is central.

So suppose $\sum_{n+m\ge 1} \alpha_{nm} C_1^nC_2^m + \sum_{n_j,m_j} C_1^{n_j}C_2^{m_j} \sum_i \beta_{ij}u_{ij}$ is central. The computation of $C_1^nC_2^m$
done just before Lemma~\ref{autoorder2} yields that the $(1,2)$ entry of $\sum_{n+m\ge 1} \alpha_{nm}C_1^nC_2^m$ will be equal to
\[
\sum_{n+m\ge 1} \alpha_{nm} (q_{n-1} x_2'{}^my_1 + Q_{m-1} x_1^ny_2) + \mu, \qquad \mu\in K[X;Y]^{(3)}.
\]
Analogously the $(1,2)$ entry of $\sum_{n_j,m_j} C_1^{n_j} C_2^{m_j} \sum_i \beta_{ij}u_{ij}$ is
\[
\sum_{n_j,m_j,i} \beta_{ij} x_1^{n_j} x_2^{m_j} (x_1'-x_1)^{r_{ij}-1} (x_2'-x_2)^{s_{ij}-1}  (y_1(x_2'-x_2) - y_2(x_1'-x_1)) + \rho.
\]
Here $r_{ij} = \deg_{t_1} u_{ij}$, $s_{ij} = \deg_{t_2} u_{ij}$, and $\rho\in K[X;Y]^{(3)}$. Since our element is central its $(1,2)$ entry lies in
$K[X;Y]^{(3)}$. Thus we obtain that the sum
\begin{eqnarray*}
&&\sum_{n+m\ge 1} \alpha_{nm}(q_{n-1} x_2'{}^m y_1+ Q_{m-1} x_1^ny_2) + \\
&&\sum_{n_j,m_j,i} \beta_{ij} x_1^{n_j}x_2^{m_j} (x_1'-x_1)^{r_{ij}-1} (x_2'-x_2)^{s_{ij}-1} (y_1(x_2'-x_2)-y_2(x_1'-x_1))
\end{eqnarray*}
must vanish. But the set $B_1$ is a basis of $K[X;Y]^{(1)}$ therefore
\begin{eqnarray*}
\sum_{n+m\ge 1} \alpha_{nm} q_{n-1}x_2'{}^m + \sum_{n_j,m_j,i} \beta_{ij} x_1^{n_j} x_2^{m_j} (x_1'-x_1)^{r_{ij}-1} (x_2'-x_2)^{s_{ij}} &=&0\\
\sum_{n+m\ge 1} \alpha_{nm} Q_{m-1} x_1^n - \sum_{n_j,m_j,i} \beta_{ij} x_1^{n_j} x_2^{m_j} (x_1'-x_1)^{r_{ij}} (x_2'-x_2)^{s_{ij}-1} &=&0.
\end{eqnarray*}
Multiplying the first equation by $(x_1'-x_1)$, the second by $(x_2'-x_2)$ and summing up we will obtain that
\begin{eqnarray*}
0&=&\sum_{n+m\ge 1} \alpha_{nm} (q_{n-1} x_2'^{m} (x_1'-x_1) + Q_{m-1} x_1^n (x_2'-x_2))\\
&=&\sum_{n+m\ge 1} \alpha_{nm}(x_2'{}^m (x_1'{}^n-x_1^n) + x_1^n(x_2'{}^m-x_2^m)) \\
&=& \sum_{n+m\ge 1} \alpha_{nm} (x_2'{}^m x_1'{}^n - x_1^nx_2^m).
\end{eqnarray*}
Therefore $\alpha_{nm}=0$ whenever $n+m\ge 1$. So we are left with the sum
\[
\sum_{n_j,m_j,i} \beta_{ij} x_1^{n_j} x_2^{m_j} (x_1'-x_1)^{r_{ij}-1} (x_2'-x_2)^{s_{ij}-1} (y_1(x_2'-x_2) - y_2(x_1'-x_1))=0.
\]
Thus we have that $\sum_{n_j,m_j,i} \beta_{ij} x_1^{n_j} x_2^{m_j} (x_1'-x_1)^{r_{ij}-1} (x_2'-x_2)^{s_{ij}} =0$. Similarly
$\sum_{n_j,m_j,i} \beta_{ij} x_1^{n_j} x_2^{m_j} (x_1'-x_1)^{r_{ij}} (x_2'-x_2)^{s_{ij}-1} =0$.
By homogeneity we deduce that for each $j$ it holds $\sum_i \beta_{ij} (x_1'-x_1)^{r_{ij}} (x_2'-x_2)^{s_{ij}} =0$. Recalling the definition of the sets $I_{ij}$ we have $\sum_{p\in I_{ij}} \beta_{pj} =0$ and we are done.
\epr

\begin{corollary}
For the centre $Z(F)$ we have $Z(F) = K\oplus I$ where $I$ is a nilpotent ideal of $F$ (and not only of $Z(F)$).
\end{corollary}

We observe that the last Corollary, together with Theorem~\ref{descriptioncentre} gives a precise answer to the question of Berele, and that $I$ is a nilpotent ideal actually of $F$, not only of the centre.

A further remark is relevant. It is interesting to note that when one deals with 3 generators, say the generic matrices $C_1$, $C_2$, $C_3$, then the element $[C_1,C_2,[C_1,C_3]]$ is central in the generic algebra of three generators. But $C_2[C_1,C_2,[C_1,C_3]]$ is not. Therefore the analogue of the above nilpotent ideal is an ideal of the centre only.

\end{document}